\documentclass[10pt,psamsfonts]{amsart}

\usepackage{amsmath,amssymb}
\usepackage{graphicx}
\usepackage[dvips]{psfrag}

\newtheorem{theorem}{Theorem}
\newtheorem{proposition}[theorem]{Proposition}

\newtheorem{lemma}[theorem]{Lemma}

\newcommand{\Cs}{\mbox{\rm Cs}}
\newcommand{\Sn}{\mbox{\rm Sn}}

\title[On the centers of the weight--homogeneous pol. vector fields
] {On the centers of the weight--homogeneous polynomial vector
fields on the plane}
\author[Jaume Llibre and Claudio Pessoa]{}

\thanks{The first author is partially supported by a MCYT
grant MTM2008--03437--C02--01 and by a CIRIT  grant number
2005SGR00550. The second author is supported by FAPESP--Brasil
grant 06/56664--0.}

  \subjclass{Primary 34C35, 58F09; Secondary 34D30}
   \keywords{Weight--homogeneous polynomial differential systems, centers, Liapunov constants}

\begin{document}
 \maketitle

\centerline{\scshape  Jaume Llibre}
\medskip

{\footnotesize \centerline{ Universitat Aut\`{o}noma de Barcelona,
Departament de Matem\`{a}tiques} \centerline{08193 Bellaterra,
Barcelona, Catalonia, Spain }
\centerline{\email{jllibre@mat.uab.cat}}}

\medskip

\centerline{\scshape  Claudio Pessoa}
\medskip

{\footnotesize \centerline{Universidade Federal de Uberl\^andia,
Faculdade de Matem\'atica}\centerline{Campus Santa M\^onica -
Bloco 1F, Sala 1F118} \centerline{Av. Jo\~ao Naves de Avila, 2121,
38.408--100, Uberl\^andia, MG, Brasil }
\centerline{\email{pessoa@famat.ufu.br}}}

\medskip

\bigskip

\begin{quote}{\normalfont\fontsize{8}{10}\selectfont
{\bfseries Abstract.} We classify all the centers of a planar
weight--homogeneous polynomial vector field of weight degree $1$,
$2$, $3$ and $4$.
\par}
\end{quote}

\section{Introduction and statement of the main results}

\noindent In the qualitative theory of real planar polynomial
differential systems two of the main problems are the
determination of limit cycles and the center--focus problem; i.e.
to distinguish when a singular point is either a focus or a
center. The notion of center goes back at least to Poincar\'e in
\cite{HP}. He defined it for a vector field on the real plane.
This paper deals with the classification of the centers for a
class of polynomial differential systems. The classification of
the centers of the polynomial differential systems started with
the quadratic ones with the works of Dulac \cite{HD}, Kapteyn
\cite{WK1}, \cite{WK2}, Bautin \cite{NB}, Zoladek \cite{HZ}, ...
see Schlomiuk \cite{DS} for an update on the quadratic centers.
There are many partial results for the centers of polynomial
differential systems of degree larger than $2$, but we are very
far to obtain a complete classification of all centers for the
polynomial differential systems of degree $\geq 3$.

We deal with polynomial differential systems of the form
\begin{equation}
\label{eq:01} \begin{array}{c} \dot{x} = P(x,y),\\ \dot{y} =
Q(x,y),
\end{array}
\end{equation}
where $P$ and $Q$ are polynomials in the variables $x$ and $y$
with real coefficients.

We say that system \eqref{eq:01} is {\it weight--homogeneous} if
there exist $s=(s_1,s_2)\in \mathbb N^2$ and $d\in \mathbb N$ such
that for arbitrary $\lambda\in \mathbb R^+=\{\lambda\in\mathbb R:
\lambda>0\}$,
\[
P(\lambda^{s_1}x,\lambda^{s_2}y)=\lambda^{s_1-1+d}P(x,y),\qquad
Q(\lambda^{s_1}x,\lambda^{s_2}y)=\lambda^{s_2-1+d}Q(x,y).
\]
We call $s=(s_1,s_2)$ the {\it weight exponent} of system
\eqref{eq:01} and $d$ the {\it weight degree} with respect to the
weight exponent $s$. In the particular case that $s=(1,1)$ systems
\eqref{eq:01} are exactly the {\it homogeneous polynomial
differential systems of degree} $d$.

A singular point $p$ of system \eqref{eq:01} is a {\it center} if
there is a neighborhood of $p$ fulfilled of periodic orbits with
the unique exception of $p$. The {\it period annulus} of a center
is the region fulfilled by all the periodic orbits surrounding the
center. We say that a center located at the origin is {\it global}
if its period annulus is $\mathbb R^2\setminus \{(0,0)\}$.

%
Let $U$ be an open subset of $\mathbb R^{2}$. Here a nonconstant
analytic function $H:U\rightarrow \mathbb R$ is called a {\it
first integral} of system \eqref{eq:01} on $U$ if it is constant
on all solutions curves $(x(t), y(t))$ of the vector field $X$
associated to system \eqref{eq:01} on $U$; i.e. $H(x(t), y(t))=$
constant for all values of $t$ for which the solution $(x(t),
y(t))$ is defined in $U$. Clearly $H$ is a first integral of the
vector field $X$ on $U$ if and only if
\[
XH=P\frac{\partial H}{\partial x}+Q\frac{\partial H}{\partial
y}\equiv 0
\]
on $U$.

The main goal of this paper is to classify all centers of the
weight--homogeneous planar polynomial differential systems
\eqref{eq:01} of weight degree $1$, $2$, $3$ and $4$ with $P$ and
$Q$ coprime.

{For} doing that we characterize first the normal forms of all the
weight--homogeneous planar polynomial differential systems of
weight degree $1$, $2$, $3$ and $4$. Thus from the definition of
weight--homogeneous polynomial differential systems \eqref{eq:01}
with weight degree $1$, $2$, $3$ and $4$, the exponents u and v of
any monomial $x^uy^v$ of $P$ and $Q$ are such that they satisfy
respectively the relations
\[
s_1u + s_2v = s_1,\qquad s_1u + s_2v = s_2,
\]
for weight degree $1$;
\[
s_1u + s_2v = s_1+1,\qquad s_1u + s_2v = s_2+1,
\]
for weight degree $2$;
\[
s_1u + s_2v = s_1+2,\qquad s_1u + s_2v = s_2+2,
\]
for weight degree $3$; and
\[
s_1u + s_2v = s_1+3,\qquad s_1u + s_2v = s_2+3,
\]
for weight degree $4$. Moreover taking into account that we only
consider the cases with $P$ and $Q$ coprime, it is easy to check
that the systems:

{\bf (i)} with weight degree $1$ are the following ones with their
corresponding values of $s$:
\begin{eqnarray}\label{eq:15} s=(1,1): & & \begin{array}{lcl} \dot{x} & = & a_{10}x+a_{01}y,\\
\dot{y} & = & b_{10}x+b_{01}y,\end{array}
\\
\label{eq:16} s=(1,p): &  & \begin{array}{lcl} \dot{x} & = & a_{10}x,\\
\dot{y} & = & b_{p0}x^p+b_{01}y,\end{array} \mbox{with}\;p\in
\mathbb N,\; p>1,
\end{eqnarray}

{\bf (ii)} with weight degree $2$ are the following ones with
their corresponding values of $s$:
\begin{eqnarray}
\label{eq:17}  s=(1,1): & & \begin{array}{lcl} \dot{x} & = & a_{20}x^2+a_{11}xy+a_{02}y^2,\\
\dot{y} & = & b_{20}x^2+b_{11}xy+b_{02}y^2,\end{array}
\\
\label{eq:18} s=(1,2): & & \begin{array}{lcl} \dot{x} & = & a_{20}x^2+a_{01}y,\\
\dot{y} & = & b_{30}x^3+b_{11}xy,\end{array}
\\
\label{eq:19} s=(2,3): & & \begin{array}{lcl} \dot{x} & = & a_{01}y,\\
\dot{y} & = & b_{20}x^2,\end{array}
\end{eqnarray}

{\bf (iii)} with weight degree $3$ are the systems with weight
degree $2$ and additionally the following ones with their
corresponding values of $s$ (see \cite{ref:5}):
\begin{eqnarray}
\label{eq:20} s=(1,1): & & \begin{array}{lcl} \dot{x} & = & a_{30}x^3+a_{21}x^2y+a_{12}xy^2+a_{03}y^3,\\
\dot{y} & = &
b_{30}x^3+b_{21}x^2y+b_{12}xy^2+b_{03}y^3,\end{array}
\\
\label{eq:21} s=(1,2): & & \begin{array}{lcl} \dot{x} & = & a_{30}x^3+a_{11}xy,\\
\dot{y} & = & b_{40}x^4+b_{21}x^2y+b_{02}y^2,\end{array}
\\
\label{eq:22} s=(1,3): & & \begin{array}{lcl} \dot{x} & = & a_{30}x^3+a_{01}y,\\
\dot{y} & = & b_{50}x^5+b_{21}x^2y,\end{array}
\end{eqnarray}

{\bf (iv)} with weight degree $4$ are the systems \eqref{eq:17},
\eqref{eq:18} with weight degree $2$ and additionally the
following ones with their corresponding values of $s$:
\begin{eqnarray}
\label{eq:36} s=(1,1): & & \begin{array}{lcl} \dot{x} & = & a_{40}x^4+a_{31}x^3y+a_{22}x^2y^2+a_{13}xy^3+a_{04}y^4,\\
\dot{y} & = &
b_{40}x^4+b_{31}x^3y+b_{22}x^2y^2+b_{13}xy^3+b_{04}y^4,\end{array}
\\
\label{eq:37} s=(1,2): & & \begin{array}{lcl} \dot{x} & = & a_{40}x^4+a_{21}x^2y+a_{02}y^2,\\
\dot{y} & = & b_{50}x^5+b_{31}x^3y+b_{12}xy^2,\end{array}
\\
\label{eq:38} s=(1,3): & & \begin{array}{lcl} \dot{x} & = & a_{40}x^4+a_{11}xy.\\
\dot{y} & = & b_{60}x^6+b_{31}x^3y+b_{02}y^2,\end{array}
\\
\label{eq:39} s=(1,4): & & \begin{array}{lcl} \dot{x} & = & a_{40}x^4+a_{01}y,\\
\dot{y} & = & b_{70}x^7+b_{31}x^3y,\end{array}
\\
\label{eq:40} s=(2,3): & & \begin{array}{lcl} \dot{x} & = & a_{11}xy,\\
\dot{y} & = & b_{30}x^3+b_{02}y^2,\end{array}
\\
\label{eq:41} s=(2,5): & & \begin{array}{lcl} \dot{x} & = & a_{01}y,\\
\dot{y} & = & b_{40}x^4,\end{array}
\end{eqnarray}

{For} studing the centers of systems (i), (ii), (iii) and (iv), we
first simplify the weight--homogeneous planar polynomial
differential systems of weight degree $3$ with weight exponent
$(1, 1)$ having $8$ parameters to some normal forms with at most
$4$ independent parameters. After using these normal forms we
characterize which of these systems have a center. The proposition
that gives us this normal forms is the following (see
\cite{ref:11} for the proof).

\begin{proposition}
\label{pro:04} For any cubic homogeneous system there exists some
linear transformation and a rescaling of the independent variable
which transforms the system into one and only one of the following
canonical forms:
\begin{eqnarray}
\label{eq:28} &&
\begin{array}{cc}
\begin{array}{lcl} \dot{x} & = & p_1x^3+p_2x^2y+p_3xy^2+\mu y^3,\\
\dot{y} & = & \mu x^3+p_1x^2y+p_2xy^2+p_3y^3,\end{array} &
\mbox{with $\mu\neq 0$;}
\end{array}
\\
\label{eq:29} &&
\begin{array}{cc}
\begin{array}{lcl} \dot{x} & = & p_1x^3+(p_2+3\alpha)x^2y+(p_3-6\alpha)xy^2-6\alpha y^3,\\
\dot{y} & = &
p_1x^2y+(p_2-3\alpha)xy^2+(p_3+6\alpha)y^3,\end{array} &
\mbox{with $\alpha=\pm 1$;}
\end{array}
\\
\label{eq:30} &&
\begin{array}{cc}
\begin{array}{lcl} \dot{x} & = & p_1x^3+p_2x^2y+(p_3+2)xy^2-4y^3,\\
\dot{y} & = &  p_1x^2y+p_2xy^2+(p_3-2)y^3,\end{array} &
\end{array}
\end{eqnarray}
\begin{eqnarray}
\label{eq:31} &&
\begin{array}{cc}
\begin{array}{lcl} \dot{x} & = & p_1x^3+(p_2-3\alpha)x^2y+p_3xy^2-6y^3,\\
\dot{y} & = & p_1x^2y+(p_2+3\alpha)xy^2+p_3y^3,\end{array} &
\mbox{with $\alpha=\pm 1$;}
\end{array}
\\
\label{eq:32} &&
\begin{array}{cc}
\begin{array}{lcl} \dot{x} & = & p_1x^3+p_2x^2y+p_3xy^2-\alpha y^3,\\
\dot{y} & = & p_1x^2y+p_2xy^2+p_3y^3,\end{array} & \mbox{with
$\alpha=\pm 1$;}
\end{array}
\\
\label{eq:33} &&
\begin{array}{cc}
\begin{array}{lcll} \dot{x} & = & p_1x^3+(p_2-3\alpha\mu)x^2y+p_3xy^2-\alpha y^3, & \mbox{with
$\alpha=\pm 1$, $\mu>-\frac{1}{3}$} \\
\dot{y} & = & \alpha x^3+p_1x^2y+(p_2+3\alpha\mu)xy^2+p_3y^3, &
\mbox{and $\mu\neq \frac{1}{3}$;}
\end{array} & \end{array}
\\
\label{eq:34} &&
\begin{array}{cc}
\begin{array}{lcl} \dot{x} & = & p_1x^3+(p_2-\alpha)x^2y+p_3xy^2-\alpha y^3,\\
\dot{y} & = & \alpha
x^3+p_1x^2y+(p_2+\alpha)xy^2+p_3y^3,\end{array} & \mbox{with
$\alpha=\pm 1$;}
\end{array}
\\
\label{eq:35} &&
\begin{array}{cc}
\begin{array}{lcl} \dot{x} & = & x(p_1x^2+p_2xy+p_3y^2),\\
\dot{y} & = & y(p_1x^2y+p_2xy+p_3y^2).\end{array} & \end{array}
\end{eqnarray}
\end{proposition}
In what follows instead of working with system \eqref{eq:20} we
shall work with the equivalent normal form
\eqref{eq:28}--\eqref{eq:35} that have atmost four parameters.

The main results of this paper is the following one.

\begin{theorem}The following statements hold.
\begin{itemize}
\item[(a)]\label{the:02} Systems \eqref{eq:15} has a center at the
origin if and only if $b_{01}=-a_{10}$ and
$a_{10}b_{01}-a_{01}b_{10}>0$. System \eqref{eq:16} has no
centers.

\item[(b)] \label{the:03} Systems \eqref{eq:17} and \eqref{eq:19}
have no centers. System \eqref{eq:18} has a center at the origin
if and only if $a_{01}b_{30}\neq 0$ and
$(2a_{20}+b_{11})^2+8(b_{30}a_{01}-b_{11}a_{20})<0$.

\item[(c)] \label{the:04} Systems \eqref{eq:21}, \eqref{eq:28},
\eqref{eq:29}, \eqref{eq:30}, \eqref{eq:31}, \eqref{eq:32} and
\eqref{eq:35} have no centers. System \eqref{eq:22} has a center
at the origin if and only if $a_{01}b_{50}\neq 0$,
$(3a_{30}+b_{21})^2+12(b_{50}a_{01}-b_{21}a_{30})<0$ and
$3a_{30}+b_{21}=0$. Systems \eqref{eq:33} and \eqref{eq:34} have a
center if and only if $p_3=-p_1$.

\item[(d)] \label{the:05} Systems \eqref{eq:36}, \eqref{eq:37},
\eqref{eq:38}, \eqref{eq:40} and \eqref{eq:41} have no centers.
System \eqref{eq:39} has a center at the origin if and only if
$a_{01}b_{70}\neq 0$ and
$(4a_{40}+b_{31})^2+16(b_{70}a_{01}-b_{31}a_{40})<0$.
\end{itemize}
\end{theorem}

The centers of the planar weight--homogeneous polynomial vector
fields of weight degree $1$ are characterized by statement $(a)$,
of weight degree $2$ by statement $(b)$, of weight degree $3$ by
statement $(c)$ and of weight degree $4$ by statement $(d)$.

\section{Homogeneous centers}

In this section we assume that $P$ and $Q$ are coprime homogeneous
polynomials of degree $m$. We want to characterize the homogeneous
systems \eqref{eq:01} which have a center at the origin. Since P
and Q are coprime it follows that the origin is the unique (real
and finite) singular point of system \eqref{eq:01}.

System \eqref{eq:01} in polar coordinates $(\rho,\theta)$ defined
by $x = \rho \cos \theta$, $y = \rho \sin \theta$, becomes
\begin{equation}
\label{eq:02} \dot{\rho} = f(\theta)\rho^m,\; \dot{\theta} =
g(\theta)\rho^{m-1}.
\end{equation}

If the origin of system \eqref{eq:01} is a center, then from the
expression of $\dot{\theta}$ it follows that the function
$g(\theta)$ is either positive, or negative for all $\theta \in
\mathbb S^1$; otherwise if $g(\theta^*) = 0$ for some $\theta^*\in
\mathbb S^1$, then the straight line $\theta = \theta^*$ is
invariant by the flow of system \eqref{eq:01}, in contradiction
with the fact that the origin is a center. So the homogeneous
polynomial $xQ(x,y)-yP(x, y)$ of degree $m+1$ has no real factors
of degree $1$, and consequently $m$ must be odd. In the rest of
this section we only consider homogeneous systems \eqref{eq:01}
having a center at the origin. So the function $g(\theta)$ is
either positive, or negative for all $\theta\in \mathbb S^1$.
Hence system \eqref{eq:01} is equivalent to the equation
\[
\frac{d\rho}{d\theta} = \frac{f(\theta)}{ g(\theta)} \rho.
\]
Separating the variables $\rho$ and $\theta$ of this equation, and
integrating between $0$ and $\theta$ we get that any solution of
the previous differential equation is
\begin{equation}
\label{eq:03} \rho(\theta)=\left\{\begin{array}{lcr} \rho(0)\exp
\left(\displaystyle\int^{\theta}_{0}\frac{f(r)}{g(r)}dr\right) &
\mbox{if} &
f(\theta)\not\equiv 0,\\
\rho(0) & \mbox{if} & f(\theta)\equiv 0.\end{array}\right.
\end{equation}

From \eqref{eq:03} it follows immediately the next proposition.

\begin{proposition}
\label{pro:02} Suppose that P and Q are coprime homogeneous
polynomials of degree $m$. Then the origin of system \eqref{eq:01}
is a global center if and only if the polynomial $xQ(x,y)-yP(x,y)$
has no a real factors of degree $1$ (in particular $m$ is odd) and
\[
\int^{2\pi}_{0}\frac{f(r)}{g(r)}dr=0.
\]
\end{proposition}
This proposition is well known, see for instance \cite{ref:11} or
\cite{ref:2}.

\section{Monodromic singularities}

In this section we assume that $P$ and $Q$ are coprime polynomials
and let $m=\max\{\deg P, \deg Q\}$. Then we say that system
\eqref{eq:01} has {\it degree} $m$.

Consider the polynomial vector field $X=(P,Q)$ associated to
system \eqref{eq:01}. The study of the local phase portrait at the
singular points of the vector field X is a problem almost
completely solved. In most of the cases one can know which is the
behavior of the solutions in a neighborhood of a singular point.
The only case that remains open is the {\it monodromic} one. In
this case the orbits turn around the singular point. The
difficulty is in distinguishing when the orbits spiral toward or
backward the singular point (i.e. when the origin is a {\it
focus}) and when the origin is a center. This problem is know as
the {\it center--focus problem}, see \cite{ref:12} for a survey on
this problem.

Suppose that the origin is a singular point of $X$. Let $\gamma
(t)$ be an orbit of $X$ defined in a neighborhood of the origin
that tends to it when $t$ tends to $+\infty$ and such that
$\lim_{t\rightarrow +\infty}\gamma (t)/\|\gamma (t)\|\in \mathbb
S^1$, where $\mathbb S^1$ is the unit circle, with $\|\cdot\|$ the
Euclidean norm. In this case $\gamma$ is said a {\it
characteristic orbit}. We may also consider orbits tending to the
origin as $t$ tends to $-\infty$ but we can always change the sign
of the parameter $t$ in system \eqref{eq:01} and assume that $t$
tends to $+\infty$. We have that the origin is a monodromic
singular point of $X$ if there is no characteristic orbit
associated to it.

We write $P(x,y)=\sum_{i=n}^{m}P_i(x,y)$,
$Q(x,y)=\sum_{i=n}^{m}Q_i(x,y)$, where $P_i$ and $Q_i$ are the
respective homogeneous part of degree $i$ of $P$ and $Q$, with
$1\leq n\leq m$. It is possible that $P_n(x, y)$ or $Q_n(x, y)$ is
null, but both of them cannot be null. A characteristic direction
for the origin of $X$ is a root $\omega \in \mathbb S^1$ of the
homogeneous polynomial $xQ_n(x, y)- yP_n(x, y)$, which can be
written as $\omega = (\cos \theta, \sin \theta)$ with $\theta\in
[0, 2\pi)$.

It is obvious that, unless $xQ_n(x, y)- yP_n(x, y)\equiv 0$, the
number of characteristic directions for the origin of $X$ is less
than or equal to $n+1$. The following well-known result relates
characteristic orbits with characteristic directions of singular
points.

\begin{proposition}
Let $\gamma (t)$ be a characteristic orbit for the origin of $X$
associated to system \eqref{eq:01} and $\omega =
\lim_{t\rightarrow +\infty}\frac{\gamma (t)}{\|\gamma (t)\|}$.
Then $\omega$ is a characteristic direction for $X$.
\end{proposition}
This proposition is proved in \cite{ref:1}. The reciprocal is not
true, see \cite{ref:12} or \cite{GGG}.

In polar coordinates $x = \rho \cos \theta$, $y = \rho \sin
\theta$, and doing a change of time $t\mapsto \tau$ such that
$dt/d\tau=\rho^{n-1}$, system \eqref{eq:01} becomes
\begin{equation}
\label{eq:43}
\begin{array}{lcl}
\dot{\rho} & = & \rho (R_n(\theta)+\rho R_{n+1}(\theta)+\cdots
+\rho^{m-n}R_m(\rho)), \\
\dot{\theta} & = & F_n(\theta)+\rho F_{n+1}(\theta)+\cdots
+\rho^{m-n}F_m(\rho)),
\end{array}
\end{equation}
where
\[
\begin{array}{c}
R_j(\theta)=P_j(\cos\theta,\sin\theta)\cos\theta+Q_j(\cos\theta,\sin\theta)\sin\theta,
\\
F_j(\theta)=Q_j(\cos\theta,\sin\theta)\cos\theta-P_j(\cos\theta,\sin\theta)\sin\theta,
\end{array}
\]
are homogeneous trigonometric polynomials of degree $j + 1$, for
$n \leq j\leq m$.

Let $\gamma (t)=(\rho (t),\theta (t))$ be an orbit for system
\eqref{eq:01} written in polar coordinates. We have that
$\gamma(t)$ tends to the origin if $\rho (t)$ is not identically
zero and $\lim_{t\rightarrow+\infty} \rho (t) = 0$. Moreover
$\gamma (t)$ tends to the origin spirally if it tends to the
origin and $\lim_{t\rightarrow+\infty} \theta (t) = \pm \infty$.
Therefore $\gamma (t)$ is a characteristic orbit for the origin of
$X$ if it tends to the origin and $\lim_{t\rightarrow+\infty}
\theta (t) =\theta^*<\infty$. In this case we say that $\theta^*$
is the tangent at the origin of $\gamma (t)$.

The polynomial $xQ_n(x, y) - yP_n(x, y)$, considered over $\mathbb
S^1$, in polar coordinates is $F_n(\theta)$ given in
\eqref{eq:43}. Hence if $\gamma (t) = (\rho (t), \theta (t))$ is a
characteristic orbit for the origin of $X$ written in polar
coordinates with tangent $\theta^*$ at the origin, then
$F_n(\theta^*) = 0$.

Now we include here a well-known property of monodromic singular
points, see \cite{GGG} for an easy proof.

\begin{proposition}
\label{pro:05} If the origin is monodromic, then
$F_n(\theta)\not\equiv 0$ and $\dot{\theta}(\rho, \theta)$ is of
definite sign for any $(\rho, \theta)$ verifying $0 < \rho <
\epsilon$, for a certain $\epsilon > 0$ sufficiently small.
\end{proposition}

When the eigenvalues of the matrix $DX(p)$ at the singular point
$p$ are complex and not real, we know that the origin is
monodromic. Such class of singular points are called of {\it
linear type}. If their real part is different from zero then the
singular point is a focus, while if their real part is zero the
singular point may be a center or a focus. This last case is the
classical Lyapunov–-Poincar\'e center problem. The study of this
problem for a concrete family of differential equations passes
through the calculation of the so--called {\it Lyapunov
constants}. The problem of how many constants are needed to
distinguish between a center of a focus is in general open.

When the matrix $DX(p)$ has its two eigenvalues equal to zero but
the matrix is not identically null, it is said that p is a {\it
nilpotent} singular point. The monodromy problem in this case was
solved in \cite{ref:13} and the center problem has been studied in
\cite{ref:14} and \cite{GGL}. Nevertheless, in general, given a
polynomial system with a nilpotent monodromic singular point it is
not an easy task to know if it is a focus or a center.

Suppose that $X$ has a singular point $p$ such that the matrix
$DX(p)$ is a nilpotent matrix. If the singularity $p$ is a center
then $p$ is called a {\it nilpotent center}. In this case using
suitable coordinates system \eqref{eq:01} can be written as
\begin{equation}
\label{eq:04}
\begin{array}{lcl}
\dot{x} & = & y + P_2(x, y),\\ \dot{y} & = & Q_2(x, y),\end{array}
\end{equation}
where $P_2$, $Q_2$ are polynomials of degree at least $2$.

We state a theorem, proved in \cite{ref:13} (see also \cite{ref:1}
and \cite{ref:15}), which solves the monodromic problem for
nilpotent singular points.

\begin{theorem}
\label{the:01} Consider system \eqref{eq:04} and assume that the
origin is an isolated singularity. Define the functions
\[
\begin{array}{lclr}
f(x) & = & Q_2(x,F(x))= ax^\alpha+O(x^{\alpha+1}), &\\
\phi(x)& = & \mbox{\rm
div}(y+P_2(x,y),Q_2(x,y))\mid_{y=F(x)}=bx^\beta+O(x^{\beta+1}), &
\end{array}
\]
where $a\neq 0$, $\alpha \geq 2$, $b\neq 0$ and $\beta \geq 1$, or
$\phi (x)\equiv 0$. Here the function $y = F(x)$ is the solution
of $y + P_2(x, y) = 0$ passing through $(0, 0)$. Then the origin
of \eqref{eq:04} is monodromic if and only if $a$ is negative,
$\alpha$ is an odd number $(\alpha = 2n-1)$, and one of the
following three conditions holds: $\beta>n-1$; $\beta=n-1$ and
$b^2+4an<0$; and $\phi\equiv 0$.
\end{theorem}

Once we know how to distinguish between monodromic and not
monodromic nilpotent singular points, we wish to solve the
center--focus problem for the monodromic ones. Our approach to
this problem passes through the computation of the Taylor
expansion of the return map near the singular point. In the next
section we define the {\it generalized Lyapunov constants}
introduced to solve the stability problem. Instead of working with
system \eqref{eq:04} we shall use the special normal form given in
the next statement, see \cite{ref:1} or \cite{ref:15} for the
proof.

\begin{lemma}
\label{le:01} An analytic vector field having an isolated
nilpotent monodromic singularity can be written as
\begin{equation}
\label{eq:05}
\begin{array}{lcl}
\dot{x} & = & y(-1+\tilde{P}(x,y)), \\ \dot{y} & = & f(x)+\phi
(x)y+\tilde{Q}(x,y)y^2,
\end{array}
\end{equation}
where all the above functions are analytic at the origin and
satisfy $\tilde{P}(0,0)=0$, $f(x)=x^{2n-1}+O(x^{2n})$ for some
$n\in \mathbb N$, $n\geq 2$, and $\phi(x)\equiv 0$ or
$\phi(x)=bx^\beta+O(x^{\beta+1})$ for some $\beta \in \mathbb N$
satisfying $\beta \geq n-1$. Moreover, if $\beta =n-1$ then
$b^2-4n<0$.
\end{lemma}

\section{Generalized Lyapunov constants}

Using the notation of previous section we will introduce some {\it
generalized polar coordinates}. As a starting point recall that
given any natural number $n\in \mathbb N$, the {\it generalized
trigonometric functions} $x(\theta) = \mbox{Cs}(\theta)$,
$y(\theta) = \mbox{Sn}(\theta)$, are defined as the unique
solution of the Cauchy problem
\begin{equation}
\label{eq:06}
\begin{array}{lcl}
\dot{x}=\displaystyle\frac{dx}{d\theta} & = & -y, \\ \\
\dot{y}=\displaystyle\frac{dy}{d\theta} & = & x^{2n-1},
\end{array}
\end{equation}
with initial conditions $x(0)=1$, $y(0)=0$.

Some properties of these functions are stated in the next
proposition, see \cite{ref:16}.

\begin{proposition}
\label{pro:01} Let $n\in \mathbb N$ and let $\Cs (\theta)$ and
$\Sn (\theta)$ be the generalized trigonometric functions
determined by system \eqref{eq:06}. Then the following
 statement hold.
\begin{itemize}
\item[(a)] $\Cs^{2n}(\theta)+n\Sn^2(\theta)=1$.

\item[(b)] $\Cs(\theta)$ and $\Sn(\theta)$ are $T$--periodic
functions where
\[
T=2\sqrt{\frac{\pi}{n}}\frac{\Gamma\left(\frac{1}{2n}\right)}{\Gamma\left(\frac{n+1}{2n}\right)}.
\]

\item[(c)] $\Sn(\theta+T/2)=-\Sn(\theta)$,
$\Cs(\theta+T/2)=-\Cs(\theta)$, $\Sn(-\theta+T/2)=\Sn(\theta)$ and
$\Cs(-\theta+T/2)=-\Cs(\theta)$.
\end{itemize}
\end{proposition}
Here as is usual $\Gamma$ denotes the {\it gamma} function, i.e.
$\Gamma(z)=\int^\infty_0e^{-t}t^{z-1}dt$.

We define for any given $n\in \mathbb N$ the {\it generalized
polar coordinates} $r$ and $\theta$ of the real plane $(x, y)\in
\mathbb R^2$ as
\begin{equation}
\label{eq:07} x=r\Cs(\theta),\; y=r^n\Sn(\theta),
\end{equation}
where the function $\Sn$ and $\Cs$ are given through
\eqref{eq:06}. Notice that $x^{2n} +ny^{2n} = r^{2n}$. Furthermore
the following equalities hold
\begin{equation}
\label{eq:08}
\begin{array}{lcl}
\dot{r} & = & \displaystyle
\frac{x^{2n-1}\dot{x}+y\dot{y}}{r^{2n-1}}, \\ \\
\dot{\theta} & = & \displaystyle
\frac{x\dot{y}-ny\dot{x}}{r^{n+1}}.
\end{array}
\end{equation}

Introducing the coordinates \eqref{eq:07}, it follows from
\eqref{eq:08} that system \eqref{eq:05} can be reduced to the
following differential equation
\begin{equation}
\label{eq:09}
\frac{dr}{d\theta}=R(r,\theta)=\frac{G(r,\theta)}{1+H(r,\theta)},
\end{equation}
where $R$ is analytic at the origin. Thus the above equation can
be written in the neighborhood of the origin as
\begin{equation}
\label{eq:10} \frac{dr}{d\theta}=\sum_{i=1}^\infty R_i(\theta)r^i,
\end{equation}
where the functions $R_i$ for $i\geq 1$ are $T$--periodic, being
$T$ the period of the generalized trigonometric functions. Now we
can define the {\it generalized Lyapunov constants} see
\cite{ref:16}.

Consider the solution $r(\theta,\rho)$ of \eqref{eq:10} such that
$r(0,\rho)=\rho$. It can be written as
\begin{equation}
\label{eq:11} r(\theta,\rho)=\sum_{i=1}^{\infty}u_i(\theta)\rho^i,
\end{equation} where $u_1(0)=1$ and $u_k(0)=0$ for all $k\geq 2$.
Hence the return map is given by the series
\[
P(\rho)=r(T,\rho)=\sum_{i=1}^{\infty}u_i(T)\rho^i.
\]
{For} a fixed system the only significant term is the first that
makes the return map differ from the identity map, and this will
determine the stability of the origin. On the other hand if we
consider a family of systems depending on parameters, each of the
$u_i(T)$ depends on these parameters. Thus the stability of the
origin is given by the first $V_k=u_k(T)\neq 0$ with $u_1(T)=1$,
$u_2(T)=\cdots=u_{k-1}(T)=0$, which it is called $k$th {\it
generalized Lyapunov constant}. Note that system \eqref{eq:04} has
a center if and only if $V_1=1$ e $V_k=0$ for all $k\geq 2$.

In \cite{ref:15} it is proved the following result.
\begin{proposition}
\label{pro:03} Consider system \eqref{eq:05} under the assumptions
of Lemma \ref{le:01}. Then the origin is a monodromic singular
point and its first generalized Lyapunov constants is
\[
V_1=\exp\left(\frac{2b\pi}{n\sqrt{4n-b^2}}\right),
\]
when $\beta = n-1$ and $n$ is odd; otherwise $V_1 = 1$.
\end{proposition}

\section{Weight--homogeneous centers of Weight degree $1$, $2$, $3$ and $4$}
In this section we shall prove the main results of this paper. We
remark that for a weight--homogeneous polynomial vector field
$X=(P,Q)$ with $P$ and $Q$ coprime the origin is the unique (real
and finite) singularity of $X$. Moreover when it is a center it is
a global center (see \cite{ref:2}).

\bigskip\noindent {\it Proof of Theorem \ref{the:02}~$(a)$.} Consider
system \eqref{eq:15}. This system is linear and we know that the
origin is a center if and only if $b_{01}=-a_{10}$ and
$a_{10}b_{01}-a_{01}b_{10}>0$. Now system \eqref{eq:16} does not
have a center because the straight line $x=0$ is invariant by its
flow.\hfill $\Box$

\bigskip\noindent {\it Proof of Theorem \ref{the:03}~$(b)$.} Consider the
vector field $X=(P,Q)$ associated to system \eqref{eq:17}. This
vector field does not have centers by Proposition \ref{pro:02},
because $m$ is even.

System \eqref{eq:19} does not have a center, because the origin is
always a cusp, by a result that can be find in \cite{ref:1} page
$362$. Note that this system has the first integral
$H(x,y)=a_{01}y^2/2-b_{20}x^3/3$.

Now we will study system \eqref{eq:18}. Note that
$a_{01}b_{30}\neq 0$, otherwise the straight lines $x=0$ or $y=0$
should be invariant by the flow of the system, and so the origin
cannot be a center. Therefore rescaling the independent variable
system \eqref{eq:18} becomes
\[
\begin{array}{lcl} \dot{x} & = & y+P_2(x,y)=\displaystyle
y+\frac{a_{20}}{a_{01}}x^2,\\\\
\dot{y} & = & Q_2(x,y)=\displaystyle
\frac{b_{30}}{a_{01}}x^3+\frac{b_{11}}{a_{01}}xy.\end{array}
\]
We have that $F(x)=-a_{20}x^2/a_{01}$ is the solution of
$y+P_2(x,y)=0$. Hence
\[
\begin{array}{lcl}
f(x) & = & \displaystyle
Q_2(x,F(x))=\frac{b_{30}a_{01}-b_{11}a_{20}}{a_{01}^2}x^3,
\\\\
\phi(x)& = & \mbox{\rm
div}(y+P_2(x,y),Q_2(x,y))\mid_{y=F(x)}=\displaystyle\frac{2a_{20}+b_{11}}{a_{01}}x.
\end{array}
\]
Thus by Theorem \ref{the:01} the origin of system \eqref{eq:18} is
monodromic if and only if
\[
(2a_{20}+b_{11})^2+8(b_{30}a_{01}-b_{11}a_{20})<0.
\]
In this case the origin is a center of \eqref{eq:18}, because the
system above is reversible, i.e. it is invariant by the change of
variables $(x,y,t)\mapsto (-x,y,-t)$.\hfill $\Box$

\bigskip\noindent {\it Proof of Theorem \ref{the:04}~$(c)$.} We begin
studying system \eqref{eq:21}. This system does not have centers,
because the straight line $x=0$ is invariant by its flow.

Consider system \eqref{eq:22}. Note that $a_{01}b_{50}\neq 0$,
otherwise the straight lines $x=0$ or $y=0$ should be invariant by
the flow of this system, and so the origin cannot be a center.
Rescaling the independent variable system \eqref{eq:22} becomes
\begin{equation}
\label{eq:26}
\begin{array}{lcl} \dot{x} & = & \displaystyle y+\frac{a_{30}}{a_{01}}x^3,\\
\dot{y} & = & =\displaystyle
\frac{b_{21}}{a_{01}}x^2y+\frac{b_{50}}{a_{01}}x^5.\end{array}
\end{equation}
In a similar way as in the study of system \eqref{eq:18} in the
proof of Theorem \ref{the:03}~$(b)$, we have that the origin of
system \eqref{eq:22} is monodromic if and only if
\[
(3a_{30}+b_{21})^2+12(b_{50}a_{01}-b_{21}a_{30})<0.
\]
Now by the change of variables
\[
\tilde{y}=\left(\frac{a_{30}b_{21}-a_{01}b_{50}}{a_{01}^2}\right)^{1/4}y+
\frac{a_{30}}{a_{01}}\left(\frac{a_{30}b_{21}-a_{01}b_{50}}{a_{01}^2}\right)^{1/4}x^3,
\]
\[
\tilde{x}=-\left(\frac{a_{30}b_{21}-a_{01}b_{50}}{a_{01}^2}\right)^{1/4}x,
\]
system \eqref{eq:26} can be written as
\begin{equation}
\label{eq:27}
\begin{array}{lcl} \dot{x} & = & \displaystyle -y,\\
\dot{y} & = & \displaystyle x^5+\mbox{sign}
(a_{01})\frac{3a_{30}+b_{21}}{\sqrt{a_{30}b_{21}-b_{50}a_{01}}}x^2y.\end{array}
\end{equation}
Therefore by Proposition \ref{pro:03} the first generalized
Lyapunov constant $V_1$ of system \eqref{eq:27} is
\[
V_1=\exp\left(\displaystyle\frac{2\left(\mbox{sign}
(a_{01})\frac{3a_{30}+b_{21}}{\sqrt{a_{30}b_{21}-b_{50}a_{01}}}\right)\pi}{n\sqrt{4n-\left(\mbox{sign}
(a_{01})\frac{3a_{30}+b_{21}}{\sqrt{a_{30}b_{21}-b_{50}a_{01}}}\right)^2}}\right).
\]
Hence $3a_{30}+b_{21}=0$ is a necessary condition for the origin
be a center of system \eqref{eq:22}. But this condition is also
sufficient, because in this case system \eqref{eq:27} has the
first integral $H(x,y)=y^2/2+x^6/6$. Note that the origin is a
global minimum of the graphic of the $H$.

{For} finishing the proof of statement $(c)$ we must study system
\eqref{eq:20}. We consider the normal forms of system
\eqref{eq:20} given by Proposition \ref{pro:04}.

Let $X=(P,Q)$ be the vector field associated to the systems given
by Proposition \ref{pro:04}. System \eqref{eq:28} does not have a
center by Proposition \ref{pro:02}, because due to the fact that
\[
xQ(x,y)-yP(x,y)=\mu(x-y)(x+y)(x^2+y^2),
\]
it has the invariant straight lines $x\pm y=0$ though the origin.

In the same away system \eqref{eq:29} does not have a center,
because
\[
xQ(x,y)-yP(x,y)=-6\alpha y^2(x^2-2xy-y^2).
\]

Now systems \eqref{eq:30}, \eqref{eq:31}, \eqref{eq:32} and
\eqref{eq:35} do not have a center, because $y=0$ is an invariant
straight line for these systems.

{For} system \eqref{eq:34} we have that
\[
xQ(x,y)-yP(x,y)=\alpha(x^2+y^2)^2.
\]
On the other hand for this system we get that
\[
\int^{2\pi}_{o}\frac{f(\theta)}{g(\theta)}d\theta=\frac{\pi(p_1+p_3)}{\alpha},
\]
where $f$ and $g$ are given by \eqref{eq:02}. Therefore, by
Proposition \ref{pro:02} system \eqref{eq:34} has a center if and
only if $p_3=-p_1$.

Analogously for system \eqref{eq:33} we have that
\[
xQ(x,y)-yP(x,y)=\alpha(x^4+6\mu x^2y^2+y^4).
\]
Hence $xQ(x,y)-yP(x,y)=0$ if and only if
\[
x=\pm |y|\sqrt{-3\mu\pm\sqrt{9\mu^2-1}}.
\]
Thus, if $\mu>-1/3$ the polynomial $x^4+6\mu x^2y^2+y^4$ does not
have a real factor of degree $1$. On another hand for this system
we get that
\[
\int^{2\pi}_{o}\frac{f(\theta)}{g(\theta)}d\theta= \displaystyle
\frac{\sqrt{3\mu+\sqrt{9\mu^2-1}}-\sqrt{3\mu-\sqrt{9\mu^2-1}}}{\sqrt{9\mu^2-1}\alpha}\pi
(p_1+p_3),
\]
where $f$ and $g$ are given by \eqref{eq:02}. Therefore, by
Proposition \ref{pro:02} system \eqref{eq:33} has a center if and
only if $p_3=-p_1$. This finish the proof of Theorem
\ref{the:04}~$(c)$. \hfill $\Box$

\bigskip\noindent {\it Proof of Theorem \ref{the:05}~$(d)$.} Consider the vector field $X=(P,Q)$ associated to system
\eqref{eq:36}. This vector field does not have centers by
Proposition \ref{pro:02}, because $m$ is even.

Systems \eqref{eq:38} and \eqref{eq:40} does not have centers,
because the straight line $x=0$ is invariant by its flow.

System \eqref{eq:41} does not have a center, because the origin is
always a cusp, by a result that can be find in \cite{ref:1} page
$362$. Note that this system has the first integral
$H(x,y)=b_{40}x^5/5-a_{01}y^2/2$.

Consider system \eqref{eq:39}. Note that $a_{01}b_{70}\neq 0$,
otherwise the straight lines $x=0$ or $y=0$ should be invariant by
the flow of this system, and so the origin cannot be a center.
Rescaling the independent variable system \eqref{eq:39} becomes
\begin{equation}
\label{eq:42}
\begin{array}{lcl} \dot{x} & = & \displaystyle y+\frac{a_{40}}{a_{01}}x^4,\\
\dot{y} & = & \displaystyle
\frac{b_{31}}{a_{01}}x^3y+\frac{b_{70}}{a_{01}}x^7.\end{array}
\end{equation}
In a similar way as in the study of system \eqref{eq:18} in the
proof of Theorem \ref{the:03}~$(b)$, we have that the origin of
system \eqref{eq:39} is monodromic if and only if
\[
(4a_{40}+b_{31})^2+16(b_{70}a_{01}-b_{31}a_{40})<0.
\]
Now if $a_{01}b_{70}\neq 0$ and the origin of system \eqref{eq:39}
is monodromic, then it is a center because system \eqref{eq:42} is
reversible, i.e. it is invariant by the change of variables
$(x,y,t)\mapsto (-x,y,-t)$.

{For} finishing the proof we must study system \eqref{eq:37}. Note
that $a_{02}b_{50}\neq 0$, otherwise the straight lines $x=0$ or
$y=0$ should be invariant by the flow of this system, and so the
origin cannot be a center.

Now in polar coordinates system \eqref{eq:37} is given by
\eqref{eq:43}, where $n=2$ and
$F_2(\theta)=-a_{02}\sin^3(\theta)$. Hence, as
$\dot{\theta}(\rho,\theta)= -a_{02}\sin^3(\theta)+\rho
F_{3}(\theta)+\cdots$, it follows that $\dot{\theta}(\rho,\theta)$
changes of sign for $0<\rho<\epsilon$ with $\epsilon$ small
enough. Therefore, by Proposition \ref{pro:05},
the origin cannot be a monodromic singular point of system
\eqref{eq:37}. Hence system \eqref{eq:37} has no centers. \hfill
$\Box$


\addcontentsline{toc}{chapter}{Bibliografia}
\begin{thebibliography}{16}

\bibitem{ref:15} {\sc M. J. \'Alvarez and A. Gasull}, Monodromy
and stability for nilpotent critical points, {\it Int. J.
Bifurcation and Chaos} {\bf 15}, 4 (2005), 1253--1265.

\bibitem{ref:13} {\sc A. F. Andreev}, Investigation of the behavior of the
integral curves of a system of two differential equations in the
neighbourhood of a singular point, {\it AMS Transl.} {\bf 8}
(1958), 183-–207.

\bibitem{ref:1} {\sc A. A. Andronov, E. A. Leontovich, I. I. Gordon and A. G. Maier},
Qualitative theory of second--order dynamic systems, John Wiley \&
Sons, $1973$.

\bibitem{NB} {\sc N.N. Bautin}, On the number of limit cycles which appear with the
variation of coefficients from an equilibrium position of focus or
center type, {\it Mat. Sbornik} {\bf 30} (1952), 181--196, {\it
AMS Transl.} {\bf 100} (1954), 1--19.

\bibitem{ref:5} {\sc L. Cair\'o and J. Llibre}, Polynomial first integrals for weight--homogeneous planar
polynomial differential systems of weight degree $3$, {\it J.
Math. Anal. and Appl.} {\bf 331} (2007), 1284--1298.

\bibitem{ref:11} {\sc A. Cima and J. Llibre}, Algebraic and topological classification of
the homogeneous cubic vector felds in the plane, {\it J. Math.
Anal. Appl.} {\bf 147} (1990), 420--448.

\bibitem{HD} {\sc H. Dulac}, D\'etermination et integration d'une certaine classe
d'\'equations diff\'erentielle ayant par point singulier un
centre, {\it Bull. Sci. Math. S\'er. (2)} {\bf 32} (1908),
230--252.

\bibitem{GGG} {\sc I. Garc\'\i a, J. Gin\'e and M. Grau}, A necessary condition in the monodromy problem for analytic
differential equations on the plane, {\it J. Symbolic Computation}
{\bf 41} (2006), 943--958.

%
%
\bibitem{GGL} {\sc H. Giacomini, J. Gin\'e and J. Llibre},
The problem of distinguishing between a center and a focus for
nilpotent and degenerate analytic systems, {\it J. Diff. Eq.} {\bf
227} (2006) 406--426.

\bibitem{WK1} {\sc W. Kapteyn}, On the midpoints of integral curves of differential
equations of the first degree, Nederl. Akad. Wetensch. Verslag.
Afd. Natuurk. Konikl. Nederland (1911), 1446--1457 (Dutch).

\bibitem{WK2} {\sc W. Kapteyn}, New investigations on the midpoints of integrals of
differential equations of the first degree, Nederl. Akad.
Wetensch. Verslag Afd. Natuurk. 20 (1912), 1354--1365; 21, 27--33
(Dutch).

%
\bibitem{ref:2} {\sc W. Li, J. Llibre, J. Yang and Z. Zhang},
Limit cycles bifurcating from the period annulus of
quasi--homogeneous centers, to appear in {\it J. Dyn. Diff. Eq.}.

\bibitem{ref:10} {\sc J. Llibre and X. Zhang}, Polynomial first integrals for
quasi--homogeneous polynomial differential systems, {\it
Nonlinearity} {\bf 15} (2002), 1269–-1280.

\bibitem{ref:16} {\sc M.A. Lyapunov}, Stability of Motion,
{\it Mathematics in Science and Engineering} {\bf 30}, Academic
Press, NY--London, 1966.

%
\bibitem{ref:12} {\sc V. Ma\~nosa}, On the center problem for degenerate
singular points of planar vector fields, {\it Int. J. Bifurcation
and Chaos} {\bf 12} (2002), 687-–707.

%
\bibitem{ref:14} {\sc R. Moussu}, Sym\'etrie et forme normale des centres et
foyers d\'eg\'en\'er\'es, {\it Ergod. Th. Dyn. Syst.} {\bf 2}
(1982), 241-–251.

\bibitem{HP} {\sc H. Poincar\'e}, M\'emoire sur les courbes d\'efinies par les \'equations
diff\'erentielles, Oeuvreus de Henri Poincar\'e, Vol. I,
Gauthiers--Villars, Paris, 1951, 95--114.

\bibitem{DS} {\sc D. Schlomiuk}, Algebraic particular integrals, integrability and
the problem of the center, {\it Trans. Amer. Math. Soc.} {\bf 338}
(1993), 799--841.

%
\bibitem{HZ} {\sc H. Zoladek}, Quadratic systems with center and their perturbations,
{\it J. Diff. Eq.} {\bf 109} (1994), 223--273.

\end{thebibliography}
\end{document}